\documentclass[12pt,a4paper,draf]{article}
\usepackage[utf8]{inputenc}
\usepackage{amsthm, amssymb, amsmath}
\usepackage[left=3cm,
right=2cm,
top=2cm,
bottom=2cm,bindingoffset=0cm]{geometry}

\usepackage[english]{babel}
\usepackage[figurename=Fig.]{caption}
\usepackage{textcomp}
\usepackage{gensymb}
\usepackage{mathrsfs}
\usepackage{xcolor} 
\usepackage{graphicx}
\graphicspath{{figures/}}
\DeclareGraphicsExtensions{.pdf,.png,.jpg}
\usepackage{multirow}

\newtheorem{Def}{Definition}[section]
\newtheorem{Prop}[Def]{Proposition}
\newtheorem{Th}[Def]{Theorem}
\newtheorem{Lem}[Def]{Lemma}
\newtheorem{Col}[Def]{Corollary}
\newtheorem{Conjecture}[Def]{Conjecture}

\theoremstyle{definition}
\newtheorem{Ex}[Def]{Example}

\newtheorem{Note}[Def]{Remark}

\theoremstyle{definition}
\newenvironment{Proof} {\par\noindent{\bf Proof. }} {\hfill$\square$\\}

\begin{document}

\title{Hamiltonian Sets of Polygonal Paths \\ in Assembly Graphs
\thanks{$^{1}$ Department of Mathematics, Bar-Ilan University, Ramat-Gan, 5290002, Israel \newline$^{2}$ Department of Mathematics and Statistics, University of South Florida, Tampa, FL 33620
\newline$^{3}$ Tel Aviv University, Tel Aviv, 6997801, Israel
\newline$^{4}$ Ben Gurion University of the Negev, Beer-Sheva, 8410501, Israel
\newline$^{5}$ HSE University, Moscow, 101000, Russia
\newline$^{6}$ Moscow Center of Fundamental and Applied Mathematics, Moscow, 119991, Russia}}
\author{A.  Guterman$^{1}$, N. Jonoska$^{2}$,  E.  Kreines$^{3,4}$,\\ A.  Maksaev$^{5}$,  N.  Ostroukhova$^{6}$}
\date{}
 
\maketitle

\begin{abstract}
We provide 
four equivalent combinatorial conditions for a simple assembly graph (rigid vertex graph where all vertices are of degree 1 or 4) to have the largest  
number of Hamiltonian sets of polygonal paths relative its size. These conditions serve 
to prove 
the conjecture 
that such 
maximum, which is equal to $F_{2n+1}-1$, where $F_k$ denotes the $k$th Fibonacci number, is achieved only for 
special assembly graphs, 
called tangled cords.  
\end{abstract}



\graphicspath{{figures/}}
\DeclareGraphicsExtensions{.pdf,.png,.jpg}
\def\occ{{\mathfrak o}}
\def\sl{{\mathfrak s}}




\section{Introduction}

Simple assembly graphs were introduced in~\cite{RNAguided_DNA_assembly}, as a tool to describe certain process of DNA recombination in some species of ciliates. 
The spacial molecular structure  at the moment of recombination at a single molecular locus is modelled by a rigid vertex of degree $4$. Multiple recombinant processes can be observed on a single DNA nanochromosome, therefore the whole process is modelled by a so-called assembly graph, a special graph with rigid vertices, see rigorous definitions in Section \ref{sec:terminology} (we also refer to~\cite{RNAguided_DNA_assembly} for the detailed and self-contained information). An assembled gene after the recombination is modelled by a polygonal path within the assembly graph, i.e., a path that visits a vertex exactly once and makes a “turn” at every vertex (except at the initial and terminal points). A set of polygonal paths that visit every vertex exactly once corresponds to a set of genes that are formed during the recombination.  Since every vertex should be visited by a polygonal path once and only once (a recombination happens only once), a set of such polygonal paths forms a so-called Hamiltonian set of polygonal paths. 
The maximal number of such Hamiltonian sets of polygonal paths gives the maximal number of sets of genes which a given molecule with specific recombination sites can encode. The question arises: given a molecule with $n$ recombination sites, what is the maximal number of sets of resulting DNA segments that can be obtained?
This translates the question above to the problem of determining the maximal number of Hamiltonian sets of polygonal paths in an assembly graph with $n$ vertices, and further characterizing those graphs where this maximum  can be achieved. 

The notion of a simple assembly graph is closely related to several mathematical 
structures, 
in particular, to   the notion of ribbon graphs or dessins d'enfants which investigates cellularly embedded graphs on oriented surfaces, see \cite{ShV}. 
There are several papers devoted to the possible genus range of the dessins d'enfants appearing from simple assembly graphs, see \cite{DNA_rec_through_as_gr_2009, Buck2012GenusRO}. 
Moreover, topological classification and enumeration of RNA structures by genus is provided in \cite{Appl}. 
In parallel, simple assembly graphs   have been studied extensively in knot theory in various contexts, see  \cite{Buck2012GenusRO} and references therein. A simple assembly graph is related to a virtual knot diagram, while a  turn at a vertex in a polygonal path corresponds to a smoothing of a vertex in the diagram~\cite{Kauff1999}. In topological graph theory, a rigid vertex graph corresponds to a graph with a fixed rotation system of the edges at every vertex, and the polygonal paths correspond to so-called A-trails, for the detailed and self-contained exposition see~\cite{A-trail}. Within this context, the Hamiltonian sets of polygonal paths in assembly graphs correspond to Hamiltonian sets of A-trails in rigid vertex graphs with all vertices of degree~4.

An approach to  characterize graphs using words, along with the notion of word-representable graphs, was first introduced by Kitaev and Pyatkin in 2008, see \cite{Kitaev2008}. Since then, this theory has undergone significant development, as evidenced by subsequent works such as \cite{Fleischmann2024, Kitaev2017, Kitaev2011, Kitaev2015} and references therein.   Sharing these ideas, a combinatorial framework linking simple assembly graphs and double occurrence words, i.e., the words in which each letter appears exactly twice, was later proposed in \cite{DNA_rec_through_as_gr_2009}.
Developing this direction further, in 2013, the number of Hamiltonian sets of polygonal paths in a simple assembly graph was estimated to be less than or equal to
$ F_{2n+1}-1$, where $F_k$ is the $k$-th Fibonacci number~\cite[Theorem 4.1]{Four_Reg_Graphs_w_Rigid_Vert_DNA_Rec_2013}, and it was shown that the bound is tight~\cite[Corollary 4.4]{Four_Reg_Graphs_w_Rigid_Vert_DNA_Rec_2013}. The given graph example achieving the bound had a specific structure, and was named a tangled cord.  
The following conjecture was stated in \cite[Conjecture 4.5]{Four_Reg_Graphs_w_Rigid_Vert_DNA_Rec_2013}:

\begin{Conjecture} \label{conj_upper_bound_only_for_TC}
The upper bound  for maximal number of Hamiltonian sets of polygonal paths, which equals   $F_{2n+1}-1$, is achieved only by the tangled cord with $n$ rigid vertices of degree 4.
	 \end{Conjecture}

 The main result of our paper is a proof of this conjecture. In order to do this we provide a characterization of the graphs with the maximal number of  Hamiltonian sets of polygonal paths in terms of combinatorics of the subwords of the corresponding double occurrence word. Our method is based on the careful analysis of the structure of the corresponding  extremal word, which leads to the characterization of the corresponding extremal graph.
 
All necessary definitions and notions such as (simple) assembly graph, polygonal paths, double occurrence words, tangled cord etc. are presented in Section~\ref{sec:terminology}.
In Section~\ref{sec:Ham-sets} we recall known facts about number of Hamiltonian sets of polygonal paths and formulate the main results of this paper. 
Proposition~\ref{Prop_equivalent_cond_F2n+1} gives four equivalent conditions for a graph that provides the largest possible number of Hamiltonian sets of polygonal paths.
Conjecture~\ref{conj_upper_bound_only_for_TC} is proven in Section~\ref{sec:main-result}.

\section{Terminology and notation} \label{sec:terminology}

\begin{figure}[h]
	\begin{center}
		\includegraphics[scale=0.65]{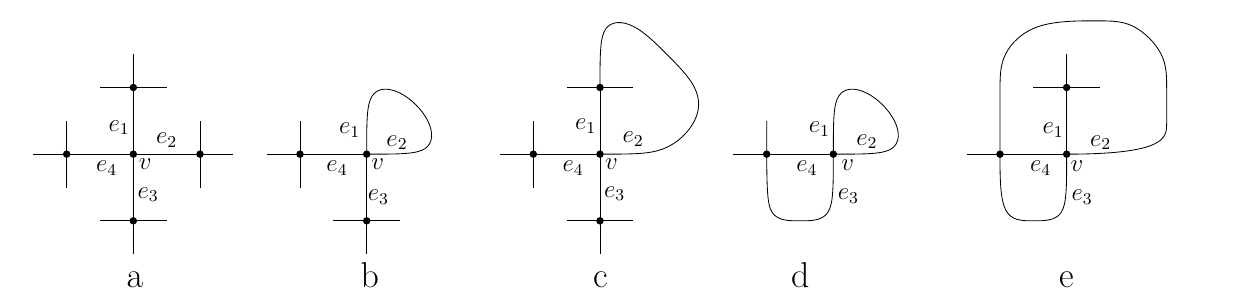}
		\caption{
			Some diagrams for a rigid vertex $v$ of degree 4.}\label{pic_ex_rigid_vertex}
	\end{center}
\end{figure}

  The definitions and notations below follow those introduced in \cite{RNAguided_DNA_assembly,DNA_rec_through_as_gr_2009, Four_Reg_Graphs_w_Rigid_Vert_DNA_Rec_2013}.  
  In this paper we consider finite undirected graphs with vertices of degree 1 or~4. Loops and multiple edges are permitted (however, multiple loops are prohibited).
  For such a graph $\Gamma$, let $V(\Gamma)$  and $E(\Gamma)$ denote the sets of vertices and edges of $\Gamma$, respectively.
 Throughout the text by degree of a vertex we mean the number of half-edges, i.e., parts of edges in a local neighborhood of a vertex, incident to this vertex, thus all multiple edges are counted  separately and a loop is counted twice.  To be more precise, if there is a loop incident to  $v \in V(\Gamma)$, the loop is represented with two half-edges (e.\,g., $e_1 $ and $ e_2$ in Figure \ref{pic_ex_rigid_vertex}b,  \ref{pic_ex_rigid_vertex}d), therefore adding two to the degree of $v$.

For a tuple of half-edges $(e_1, e_2, e_3, e_4)$, we define its {\it cyclic order} as the following set consisting of all its cyclic permutations and their reverses:
\begin{align*}
(e_1, e_2, e_3, e_4)^{cyc} = \{ & (e_1, e_2, e_3, e_4), (e_2, e_3, e_4,e_1), (e_3, e_4,e_1, e_2), (e_4,e_1, e_2, e_3), \\
& (e_4, e_3, e_2, e_1), (e_3, e_2, e_1, e_4), (e_2, e_1, e_4, e_3), (e_1,e_4,e_3,e_2)
	\}.
\end{align*}
In other words, we consider all elements of the set $(e_1, e_2, e_3, e_4)^{cyc}$ equivalent and usually refer to a single element of this ordering, mostly $(e_1, e_2, e_3, e_4)$. 
We say that a vertex $v \in V(\Gamma)$ of degree 4 is {\it rigid} if a cyclic order of  half-edges incident to this vertex is fixed (in fact there are three ways to do this). Throughout the text we depict graphs as  locally embedded to the surface at each vertex. Each vertex $v$ is considered as a small disk such that incident edges are attached to the boundary of the disk. 
The cyclic order of half-edges is usually specified as we read them on the diagram following the boundary of the disk. 
The most general type of a rigid vertex is depicted in Figure \ref{pic_ex_rigid_vertex}a,  however, there are many other possibilities, where loops and multiple edges are involved (see, e.\,g., Figure \ref{pic_ex_rigid_vertex}b-e).

Suppose that for a vertex $v$ of degree 4 the following cyclic order of half-edges is specified: $(e_1, e_2, e_3, e_4)$. Then we say that $e_1$ and $e_2$ are {\it neighbours in} $v$, and so are $e_2$ and $e_3$, and also $e_3$ and $e_4$, as well as $e_4$ and $e_1$.  Note that half-edges of a loop can be neighbours as in Figure \ref{pic_ex_rigid_vertex}b,d.

An \textit{assembly graph} is a finite connected graph in which all vertices are of degree $1$ or $4$ and all vertices of degree $4$ are rigid. A vertex of degree $1$ is called an {\it endpoint}.

\textit{A path} in a graph is an alternating sequence of vertices and edges that starts and ends with a vertex, and the edges in this sequence are incident to the vertices that stand next to them. Throughout the text, we consider that a path $v_0e_1v_1 \dots v_{\ell-1}e_{\ell}v_\ell$ and its reverse $v_\ell e_{\ell}v_{\ell-1} \dots v_1e_1v_0$ are equal.
If there are no edges in a path and so the path consists of a single vertex, it is called  \textit{a singleton}.

\begin{Def} \cite[Definition 3.3]{DNA_rec_through_as_gr_2009}
A \textit{transverse path} in an assembly graph or simply a \textit{transversal} is a path $v_0 e_1 v_1 e_2 \ldots e_n v_n$,  where $v_0, v_n$ are endpoints, satisfying the following conditions:
\begin{enumerate}
    \item All the edges $e_1, \dots, e_n$ are pairwise distinct.
    
    \item Each \(e_i\) is not a neighbour of \(e_{i-1}\) in 
    the rigid vertex \(v_{i-1}\) 
    (i = 2, \ldots, n). 
    When 
    $e_i$ is a loop ($v_{i-1} = v_i$), 
    both half-edges of $e_i$ are 
    marked as $e_i^{(1)}$ and $ e_i^{(2)}$ in such a way that $e_{i-1}$ and $e_i^{(1)}$ are not neighbours in $v_i$, as well as $e_i^{(2)}$ and $e_{i+1}$. 
\end{enumerate}
\end{Def}

\begin{figure}[ht]  
	\begin{center}
		\includegraphics[scale=1]{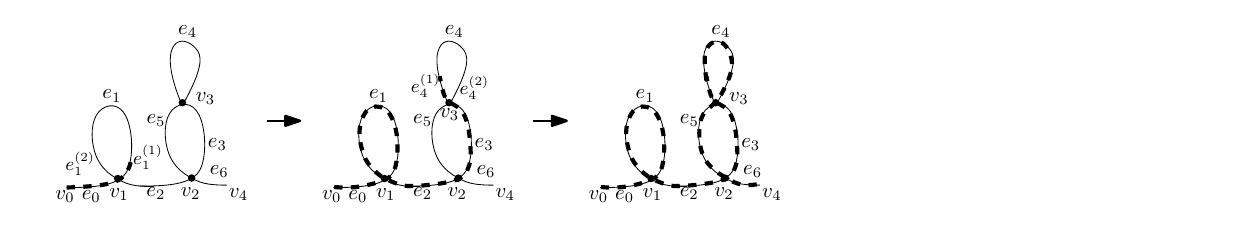}
		\caption{Consecutevily drawing transverse path $v_0\,e_0\,v_1\,e_1\,v_1\,e_2\,v_2\,e_3\,v_3\,e_4\,v_3\,e_5\,v_2\,e_6\,v_4$.}\label{pic_transversal_path}
	\end{center}
\end{figure}

An example of a transverse path can be seen in Figure \ref{pic_transversal_path}.  Note that at 
the vertex $v_1$ the transverse path first follows “right” half-edge, corresponding to the loop $e_1$, further going through the “left” half-edge of $e_1$. Meanwhile at 
$v_3$ the path first follows “left” half-edge of the loop $e_4$, further going through its “right” half-edge. 
 
Two transverse paths are \textit{equivalent} if they are either identical or one is the reverse of the other.
A \textit{simple assembly graph} is an assembly graph having a Eulerian transverse path, i.\,e., a transverse path visiting each edge exactly once. Note that in a simple assembly graph there is a unique equivalence class of transverse Eulerian paths, see \cite[Lemma~3.6]{DNA_rec_through_as_gr_2009}. Following \cite{DNA_rec_through_as_gr_2009} we say that two graphs are {\it isomorphic as simple assembly graphs} if 
there exists a graph 
isomorphism between the two graphs
that extends 
to a 
bijection on 
half-edges    preserving their respective cyclic order 
at 
all rigid vertices. In particular, a simple assembly graph is isomorphic to the simple assembly graph obtained from the original by taking all edges along the transversal in reverse order.
This means that simple assembly graphs are uniquely determined by the transverse Eulerian path, see \cite[Definition 3.2 and Lemma 3.7]{DNA_rec_through_as_gr_2009}. 

\begin{figure}[t]  
	\begin{center}
		\includegraphics[scale=0.7]{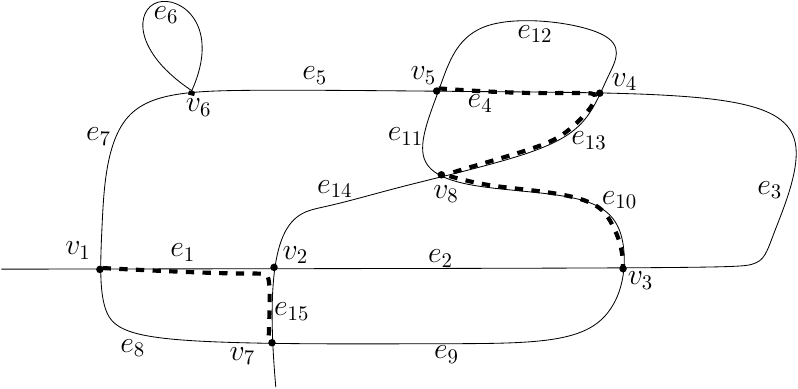}
		\caption{An example of Hamiltonian set of polygonal paths with three paths: $\{v_1 e_1 v_2 e_{15} v_7, \, v_5 e_4 v_4 e_{13} v_8 e_{10} v_3, \,v_6\}$. Note that $v_6$ is a singleton. }\label{pic_Ham_set}
	\end{center}
 \end{figure}

A path $v_0e_1v_1 \dots v_{\ell-1}e_{\ell}v_\ell$ in an assembly graph is called {\it polygonal} (see \cite[Section 2]{Four_Reg_Graphs_w_Rigid_Vert_DNA_Rec_2013}) if all $v_0, \dots, v_{\ell}$ are pairwise distinct 4-degree vertices of $\Gamma$ and $e_{i}$ and $e_{i+1}$ are neighbours in $v_{i}$, for all $i = 1, \dots, \ell-1$.  
Note that a cycle (particularly, a loop) is not a polygonal path since all the vertices in it must be  distinct, while a singleton vertex is a polygonal path. Examples of polygonal paths are presented in Figure \ref{pic_Ham_set}.

Informally, if we consider a vertex as a crossroad, we can either turn left or right forming a polygonal path or drive through at every vertex, which gives us a transverse path.

Two paths are called \textit{vertex-disjoint} if they have no common vertex.
A pairwise vertex-disjoint set of polygonal paths in an assembly graph $\Gamma$ is called \textit{Hamiltonian} (see \cite[Section 2]{Four_Reg_Graphs_w_Rigid_Vert_DNA_Rec_2013}) if its union covers all $4$-degree vertices of $\Gamma$ (see Figure~\ref{pic_Ham_set}).
We denote by $\mathcal{C}(\Gamma)$ the collection of all Hamiltonian sets of polygonal paths in $\Gamma$.

Let $\Sigma =\{a_1, a_2, \ldots\}$ be a (finite) alphabet and $w = w_1 w_2 \ldots w_{N-1} \, w_{N}$ be a word of length $N$ in the alphabet $\Sigma$. The length of $w$ is denoted by $|w|$, so $|w| = N$.
We set $w[i] = w_i$, that is, the $i$th symbol of the word $w$. 
By $w[i:j]$ $(i\le j)$ we denote the {\it subword} $w_i \, w_{i+1} \, \ldots  \, w_{j-1} \, w_j$ and
we have $|w[i:j]| = j - i +1$.
 Let $a \in \Sigma$ occur exactly twice in a fixed word $w$, in 
positions $i$ and $j$, i.\,e., $w[i] = w[j] = a$, $i < j$. We define two functions $\occ_1$ and $\occ_2$ that give the first and second occurrence of $a$ in $w$, respectively, that is, $\occ_1(a) = i$ and $\occ_2(a) =j$. In the following, the use of $\occ_1(a)$ or $\occ_2(a)$ 
presumes that the word $w$ is clear from the context and the letter $a$ appears in $w$ exactly twice.

An \textit{assembly word}, or a {\it double occurrence word} (DOW), in the alphabet $\Sigma$ is a non-empty word in which each symbol of the alphabet occurs either exactly twice or does not occur at all. Note that this combinatorial definition is used for Gauss words in the theory of plane curves. Double occurrence words are said to be {\it equivalent} if they can be obtained one from another by renaming the alphabet symbols and/or writing the word in the reverse order.
There is a one-to-one correspondence between the isomorphism classes of simple assembly graphs and the equivalence classes of double occurrence words, as proved in \cite[Lemma 3.8]{DNA_rec_through_as_gr_2009}. Below we briefly provide this construction for the completeness and further use.

Listing the vertices of degree~4 visited by a transverse path in a simple assembly graph gives a DOW over the alphabet $V(\Gamma)$, since every transverse path must visit every  vertex of degree~4 twice. 
Hence, assembly words are equivalent if and only if they correspond to isomorphic assembly graphs, as stated in  \cite[Lemma 3.8]{DNA_rec_through_as_gr_2009}.

 For a simple assembly graph $\Gamma$ we denote a representative of the corresponding class of  assembly words by $w_{\Gamma}$. Vice versa, for an assembly word $w$ we denote the corresponding simple assembly graph by  $\Gamma_w$, see Figure \ref{pic_inc_mat_constr}. 
 
 The number of degree 4 vertices in a simple assembly  graph $\Gamma$ is denoted by $\mathcal{V}(\Gamma)$. We have $ 2\cdot \mathcal{V}(\Gamma_w) = |w|$ since $w$ is a DOW. 
 We also denote cardinality of a set $X$ by $|X|$.
 
\begin{Def}\label{delete_symbols}
Let $w$ be  a word in the alphabet $\Sigma$ 
and $\sigma \subseteq \Sigma$  be a non-empty subset of letters. Denote by $w\setminus \sigma$  the sequence of non-empty subwords  obtained from $w$ by deleting the letters belonging to $\sigma$ (notice that those subwords are not necessarily assembly words even if $w$ is assembly).  
Denote by $w(\sigma)$ the concatenation of all occurrences of letters from $\sigma$ in $w$, keeping the order in which the letters occur in $w$.
 \end{Def}

	 \begin{Ex}
	 	For the word 
   $w = 134 2 134 85 67 5 7 28 6$  and $\sigma = \{2,5,8\}$, we have:
	 	$$w\setminus \sigma = (134,\, 134, \, 67,\, 7, \, 6), \, w(\sigma) = 285528.$$
	 \end{Ex}
  
\begin{figure}[h]
	\begin{center}
		\includegraphics[scale=1]{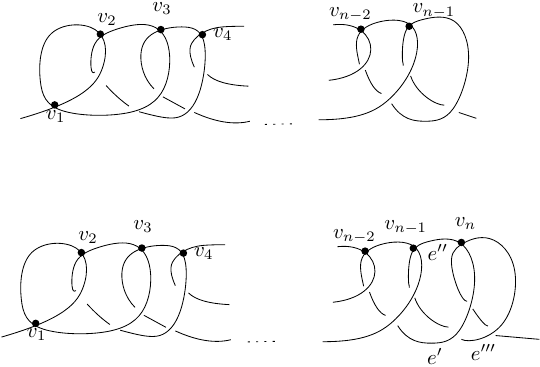}
		\caption{How the tangled cord  $TC_{n}$ is inductively constructed from $TC_{n-1}$. The visualization provides the 
        `correct' cyclic order of half-edges in each rigid vertex.}  \label{pic_add_vert_TCn}
	\end{center}
 \end{figure}

\begin{Def}\textup{\cite[Definition 4.2]{Four_Reg_Graphs_w_Rigid_Vert_DNA_Rec_2013}}\label{def_TC}
	 	A \emph{tangled cord} of order $n$ is the simple assembly graph that corresponds to the following double occurrence word:
	 	\[  w_{TC_n}= 1213243\dots (n-1)(n-2)n(n-1)n.
	 	\]
	 	The tangled cord is denoted by $TC_n$.
	 \end{Def}	
  	 Each next word $ w_{TC_n} $ is obtained from the previous assembly word $ w_{TC_{n-1}} $ via replacing the last occurrence of the letter $ (n-1) $ with a subword $ n(n-1)n$. This operation in terms of assembly graphs is represented in Figure~\ref{pic_add_vert_TCn}.

	 These are examples of assembly words $w_{TC_n}$ for small $n$: 
  $$ w_{TC_1}=11,\ w_{TC_2}=1212,\ w_{TC_3}=121323, \ w_{TC_4}=12132434.$$

\begin{Note}\label{note_tc_n_is_symmetric}
Note that the word $w_{TC_n}$ is symmetric, i.e., if we write it in reverse order and rename symbols in ascending order we obtain the same word.
\end{Note}
  
	\section{Hamiltonian sets of polygonal paths}
 \label{sec:Ham-sets}

We denote the Fibonacci sequence by $(F_n)$, namely, it is the sequence defined by $F_n = F_{n-1}+F_{n-2}$ for $n\ge 2$ and  $F_0=0, F_1=1$.
The following upper bound is known for the number of Hamiltonian sets of polygonal paths in a simple assembly graph.  
	
	\begin{Th}\textup{\cite[Theorem 4.1]{Four_Reg_Graphs_w_Rigid_Vert_DNA_Rec_2013}}\label{th_hpp_upper_bound} 
 Let $\Gamma$ be a simple assembly graph, $\mathcal{V}(\Gamma)=n $, and 
		$\mathcal{C}$ be a collection of Hamiltonian sets of polygonal paths in $\Gamma$. Then 
		$|\mathcal{C}| \leqslant F_{2n+1}-1.$
	\end{Th}

	 It was   proved in \cite[Theorem 4.3, Corollary 4.4]{Four_Reg_Graphs_w_Rigid_Vert_DNA_Rec_2013} that for the tangled cord $TC_n$ the equality $|\mathcal{C}| = F_{2n+1}-1$ is true, and therefore, the bound is tight for all $n$. It was also conjectured in the same paper \cite[Conjecture~4.5]{Four_Reg_Graphs_w_Rigid_Vert_DNA_Rec_2013} that the tangled cord is the only graph, for which  the equality holds, see  Conjecture \ref{conj_upper_bound_only_for_TC}.

In the definition below we follow the notation introduced in \cite[Remark 2.2]{DOW_Graphs_Matr}. Namely, we enumerate all the edges in the same order as they occur when going through the transversal. 

\begin{Def} \label{def:enumeration}
Let $\Gamma$ be a simple assembly graph, $\mathcal{V}(\Gamma) = n$, and let us fix the word $w_\Gamma$ by listing all the vertices of degree 4 along the 
Eulerian transversal of $\Gamma$. 
Then we enumerate the edges of $\Gamma$ by the following rule: $e_i = (w_\Gamma[i], w_\Gamma[i+1]), i = 1,2, \dots, 2n~-~1.$ 
Note that the edges incident to the endpoints of $\Gamma$ are not enumerated.
\end{Def}

\begin{Ex}
Figure \ref{pic_inc_mat_constr} shows a graph with the assembly word $w= 112323 $ and edges labelled in the order they are encountered 
along the transversal.
\begin{figure}[h] 
	\begin{center} 
		\includegraphics[scale=1.3]{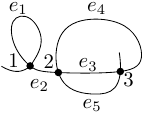} 
	\end{center} 
	\caption{Enumerated edges in a graph with the assembly word $ 112323 $
	} 
	\label{pic_inc_mat_constr} 
\end{figure}
\end{Ex}

We define a correspondence between Hamiltonian sets of polygonal paths of a simple assembly graph $\Gamma$ and  binary strings (words in the alphabet~$\{0, 1\}$).  

\begin{Def}\label{def_phi_binary_strings}
    Let $\Gamma$ be a simple assembly graph, $\mathcal{V}(\Gamma) = n$, and $e_1, e_2, \dots,$ $e_{2n-1}$ be the enumerated edges of its Eulerian transversal, see Definition \ref{def:enumeration}.
    We consider the map $ \Phi_\Gamma\colon \mathcal{C}(\Gamma)\rightarrow \{0,1\}^{2n-1}$ such that for  $\gamma \in \mathcal{C}(\Gamma)$ and $i = 1, 2, \dots, 2n-1$,  
    $$\Phi_\Gamma(\gamma)[i] = 
        \begin{cases}
            1, & \text{if $e_i$ belongs to some path from $\gamma$;}\\
            0, & \text{otherwise.}
        \end{cases}$$   
\end{Def}

\begin{Prop}\label{prop_phi_properties}
    For a simple assembly graph $\Gamma$ with $\mathcal{V}(\Gamma) = n$, the map $\Phi_{\Gamma}$ satisfies the following properties:
    \begin{enumerate}
        \item \label{prop_phi_properties_injectivity} $\Phi_{\Gamma}$  is injective.
        \item \label{prop_phi_properties_no_consequent_1} $\Phi_{\Gamma}(\gamma)$ is a string with no two consecutive ones for any $\gamma \in \mathcal{C}(\Gamma)$.
        \item \label{prop_phi_properties_1010101} 
        $\Phi_\Gamma(\gamma)\neq
    \underbrace{10101\ldots01}_{2n-1}$ 
    for any $\gamma \in \mathcal{C}(\Gamma)$.
    \end{enumerate}    
\end{Prop}
\begin{Proof}
(1)~Assume the contrary. Let $\gamma, \gamma' \in \mathcal{C}(\Gamma)$ be such that $\Phi_{\Gamma}(\gamma) = \Phi_{\Gamma} (\gamma')$. Then $\gamma$ and $\gamma'$ consist of exactly the same edges. Since both sets of paths are Hamiltonian (and in particular, vertex-disjoint), we get that $\gamma$ and $\gamma'$ coincide up to singletons.
   Each vertex in $\Gamma$ that is not incident to any edge in $\gamma$ (respectively, $\gamma'$) must be a singleton in $\gamma$ (respectively, $\gamma'$) and hence $\gamma =\gamma'$.
 
(2)~Let  $\gamma = \{\gamma_1, \gamma_2, \dots, \gamma_s\}$ be a Hamiltonian set of polygonal paths in $\Gamma$.   
Consider all edges belonging to 
 paths in $\gamma$. We claim that among these edges there are no two consecutive edges in a Eulerian  transversal of $\Gamma$.  Indeed, let $u,v,w$ be vertices of $\Gamma$ such that the edges $(u, v) \in \gamma_i$ and $(v, w) \in \gamma_j$ are consecutive in the transversal. Then $i = j$, since otherwise $\gamma_i$ and $\gamma_j$ have $v$ as a common vertex, which contradicts the definition of a Hamiltonian set. 
Therefore both edges belong to the same path: $(u, v) \in \gamma_i$ and $(v, w) \in \gamma_i$, which contradicts the polygonality of $\gamma_i$ since by the assumption $(u, v)$ and $(v, w)$ are consecutive. So, there are no consecutive edges among all edges belonging to any path in $\gamma$. Therefore there is no binary string with two consecutive $1$ in the image of~$\Phi_{\Gamma}$.

(3)~Observe that a polygonal path with $k$ vertices visits $k-1$ edges. Therefore a Hamiltonian set of polygonal paths in $\Gamma$, which visits $n$ vertices, contains at most $n-1$ edges and so would correspond to a binary string with at most $n-1$ occurrences of $1$s.  The binary string 
 $\underbrace{10101\ldots01}_{2n-1}$ has $n$ occurrences of $1$s. 
\end{Proof}

For any edge $e \in E(\Gamma)$, we denote by $v_1(e),v_2(e)\in V(\Gamma)$ the vertices incident to $e$, meaning that $e = (v_1(e),v_2(e))$. In particular, if $e$ is a loop, we still write $v_1(e),v_2(e)$ for the unique vertex of $e$. 
  
 \begin{Prop}\label{Prop_equivalent_cond_F2n+1}
	 	Let $\Gamma$ be a simple assembly graph, $\mathcal{V}(\Gamma) = n$, and $e_1, e_2, \dots,$ $e_{2n-1}$ be the enumerated edges of its Eulerian transversal in the sense of Definition \ref{def:enumeration}. Let us fix the word $w_\Gamma$ in the alphabet $\{1, 2, \dots, n\}$.   
        Then the following statements are equivalent:
    \begin{itemize}
	\item[(1)] The number of all Hamiltonian sets of polygonal paths for the graph $\Gamma$ is maximal, i.e., $|\mathcal{C}(\Gamma)| = F_{2n+1}-1.$  
	\item[(2)] Any $1 \leqslant k \leqslant n-1$ pairwise non-consecutive edges $e_{i_1},\ldots,e_{i_k}$ form a vertex disjoint set of polygonal paths.
	\item[(3)]  For any $1 \leqslant k \leqslant n-1$ pairwise non-consecutive edges $e_{i_1},\ldots,e_{i_k}$, 
        there are vertices of $\Gamma$ in the list 
        $(v_1(e_{i_1}),v_2(e_{i_1}),v_1(e_{i_2}),v_2(e_{i_2}), \ldots, v_1(e_{i_k}),v_2(e_{i_k}))$ 
        that are included into the sequence 
        exactly once.
        
	\item[(4)] For any proper non-empty subset $\sigma \subset \{1,2,\dots,n\}$, at least one word from $w_\Gamma\setminus \sigma$ has an odd length.
\end{itemize}
	 \end{Prop}

	 \begin{Proof}
$(1)\Leftrightarrow (2)$ The equivalence of Item (1) and Item (2) has been proven in \cite[Theorem 4.1]{Four_Reg_Graphs_w_Rigid_Vert_DNA_Rec_2013}, however, we provide it here for the completeness.
 By Proposition \ref{prop_phi_properties}, Item \ref{prop_phi_properties_injectivity}, the map $\Phi_{\Gamma}$ is injective. Moreover, according to Proposition \ref{prop_phi_properties}, Item \ref{prop_phi_properties_no_consequent_1}, there is no binary string with two consecutive ones in the image of $\Phi_{\Gamma}$.
            Thus we can state that the number of Hamiltonian sets of polygonal paths does not exceed the number of binary strings of length $2n-1$ with no consecutive ones. 
	 	It is well known that the number of binary strings of length $m$ with no consecutive ones  equals $F_{m+2}$ (see for example \cite[Sequence A000045]{oeis}).  
By Proposition \ref{prop_phi_properties}, Item \ref{prop_phi_properties_1010101}, the binary string $\underbrace{10101\ldots01}_{2n-1}$ does not belong to the image of  $\Phi_{\Gamma}$.

Therefore, $|\mathcal{C}(\Gamma)| \leqslant F_{2n+1}-1$, and the equality holds if and only if  any $k \leqslant n-1$ pairwise non-consecutive edges of the Eulerian transversal of $\Gamma$ form a vertex-disjoint set of polygonal paths in $\Gamma$.
   
Indeed, $|\mathcal{C}(\Gamma)| = F_{2n+1}-1$ if and only if the image of $\Phi_\Gamma$ equals all the binary strings of length $2n-1$ with no consecutive ones, excluding the string $\underbrace{10101\ldots01}_{2n-1}$. 
   This means that any binary string in $\{0, 1\}^{2n-1}$ with $k \leqslant n-1$ ones and without consecutive ones has a preimage $\gamma \in \mathcal{C}(\Gamma)$ under $\Phi_\Gamma$, which is equivalent to the fact that any
     choice of $k \leqslant n-1$ non-consecutive edges from the Eulerian transversal of $\Gamma$ provides a vertex-disjoint set of polygonal paths $\gamma'$
     ($\gamma$ is obtained from $\gamma'$ by adding all non-visited vertices as singletons).
	 	
$(2)\Rightarrow (3)$ Choose $1 \leqslant k \leqslant n-1$ arbitrary pairwise non-consecutive edges from the transversal of $\Gamma$. These edges form a vertex-disjoint set of polygonal paths in $\Gamma$. Consider any path $\alpha$ in this set. By the definition of a polygonal path, the beginning of $\alpha$ cannot belong to the other paths. Hence  it occurs exactly once among all the  endpoints of $e_{i_1},\ldots,e_{i_k}$.
	 	
$(3)\Rightarrow (2)$ Suppose  that the statement of  Item (3) is true. 
	 	Choose  $1 \leqslant k \leqslant n-1$ arbitrary pairwise non-consecutive edges $e_{i_1},\ldots, e_{i_k}$ of the Eulerian transversal of $\Gamma$.
        These edges together with the vertices incident to them induce a subgraph $\Gamma'$ of $\Gamma$, possibly, disconnected. 
	 	Since the edges $e_{i_1},\ldots,e_{i_k}$ are pairwise non-consecutive, their  endpoints  have only degrees 1 or 2 in  $\Gamma'$. Therefore, $\Gamma'$  is a union of vertex-disjoint polygonal paths and cycles (without repeating vertices).
      	
	 	If $\Gamma'$ has no cycles, we  have a vertex-disjoint set of polygonal paths.
    Assume that $\Gamma'$ has a cycle $C$. Then $C$ consists of some $e_{j_1},\ldots, e_{j_l}\in \{e_{i_1},\ldots, e_{i_k}\}$ and all their endpoints. Hence for the edges $e_{j_1},\ldots, e_{j_l}$ it holds that $1\le l\le k\le n-1$, but the sequence $(v_1(e_{j_1}),v_2(e_{j_1}),   \ldots, v_1(e_{j_l}),v_2(e_{j_l}))$ contains each vertex twice (since $C$ is a cycle), which contradicts the conditions of  Item~(3).
		 	
$(3)\Rightarrow (4)$ 
Suppose  by contradiction that all of the words from $w_\Gamma \setminus \sigma$ have even length. From each word $u = u_1u_2\ldots u_{2t} \in w_\Gamma \setminus \sigma$ we choose $t$ non-consecutive edges of the Eulerian transversal of $\Gamma$, namely, $ (u_1,u_2)$, $ (u_3,u_4),\  \dots, \ (u_{2t-1}, u_{2t})$. Since $|\sigma| \geqslant 1$, we have no more than $(2n-2)/2 = n-1$ edges chosen in total. Also, every letter from $\{1, 2, \dots, n\} \setminus \sigma$ occurs exactly twice among all endpoints of the chosen edges, contradicting Item~(3).

$(4)\Rightarrow (3)$ 
Assume, by contradiction,  that Item (3) does not hold, i.\,e., there exists $ k$, $1 \leqslant k \leqslant n-1$, and there exists a set $\mathcal{X} =\{e_{i_1}, e_{i_2}, \dots, e_{i_k}\}\subset E(\Gamma)$  such that the edges from $\mathcal{X}$ are pairwise non-consecutive in the transversal and each element in the sequence $$ \mathcal{Y}=(v_1(e_{i_1}),v_2(e_{i_1}),v_1(e_{i_2}),v_2(e_{i_2}), \ldots, v_1(e_{i_k}),v_2(e_{i_k}))$$ occurs exactly twice. 
   
   To obtain a contradiction we construct a proper non-empty subset $\sigma \subset \{1,2,\dots,n\}$ such that all words from $w_\Gamma \setminus \sigma$ have even length. 
   
   We denote  by $\delta$ the set of all distinct elements from $\mathcal{Y}$. Let us consider the subset $\sigma = \{1,2,\dots, n\} \setminus \delta$, i.\,e., $\sigma$ is the set of all 4-degree vertices of $\Gamma$ that are not endpoints of the edges from~$\mathcal{X}$. Note that $|\sigma| = n - k$.
   
   Consider an arbitrary word $u \in w_\Gamma \setminus \sigma$. We 
   show that $|u|$ is even. Let $u$ be represented as a sequence of letters $ u= u_1u_2\ldots u_{\ell}$. Note that all the letters $u_1, u_2, \dots, u_{\ell}$ correspond to vertices that belong to edges from $\mathcal{X}$, moreover, from the edges $(u_{j-1},u_j)$ and $(u_{j},u_{j+1})$ exactly one belongs to $\mathcal{X}$, for every $j = 2, \dots, \ell - 1$.
   
   The word $u$ starts with $u_1$, hence $(u_1,u_2) \in \mathcal{X}$. However, the word $u$ is a subword of  $w_\Gamma$, and by assumption the edges in $ \mathcal{X}$ are non-consecutive in the Eulerian transversal of $\Gamma$, therefore, $(u_2,u_3) \notin \mathcal{X}$. Further, $u_3$ belongs to some edge from $\mathcal{X}$, therefore $(u_3,u_4) \in \mathcal{X}$. Continuing this argument, we have the alternation of the inclusions:
	 	\[
	 	(u_1, u_2) \in \mathcal{X},\ (u_2,u_3) \notin \mathcal{X}, \ (u_3,u_4) \in \mathcal{X} \dots
	 	\] 
	 	The symmetric argument provides:
	 	\[
	 	(u_\ell, u_{\ell - 1}) \in \mathcal{X},\ (u_{\ell - 1}, u_{\ell - 2}) \notin \mathcal{X}, \ (u_{\ell - 2},u_{\ell - 3}) \in \mathcal{X} \dots 
	 	\] 
	 	Therefore, $\ell = |u|$ is even. Since $u$ was an arbitrary word from $w_\Gamma \setminus \sigma$, we obtain a contradiction. This completes the proof.
	 \end{Proof}

\begin{Def} \label{def_max_word}
    An assembly word $w$ is \emph{maximal} if $|\mathcal{C}(\Gamma_{w})| = F_{2n+1}-1$.
\end{Def}

According to Proposition \ref{Prop_equivalent_cond_F2n+1}, Item (4), an assembly word is maximal if and only if 
after deleting any non-empty proper subset of symbols from the word the remaining sequence of subwords contains at least one subword of an odd length.

\begin{Def}
	A \emph{composition} of two assembly words $u$  and $v$ without common letters (denoted by $u \circ v$) is simply a concatenation $ w = uv$.
\end{Def}

\begin{Note}
A composition of two assembly words $u$ and $v$ can be interpreted in terms of simple assembly graphs as follows (see also \cite[Definition 3.9]{DNA_rec_through_as_gr_2009}). We choose an orientation of the Eulerian transversal for each graph induced by its assembly word. Now, a composition of two directed simple assembly graphs $\Gamma_u$ and $\Gamma_v$  is the directed simple
assembly graph obtained by identifying the terminal vertex of $\Gamma_u $ with the initial vertex of $\Gamma_v$. This composition is denoted by $\Gamma_u\circ \Gamma_v$. 
\end{Note}

The following is an observation that the
composition of two assembly words cannot be maximal.

\begin{Col}\label{Col_max_not_composition}
	If $w = u\circ v$, where $u$ and $v$  are both assembly words, then $w$ is not maximal. 
\end{Col}
\begin{Proof}
	Suppose $w$ is maximal. Let $\sigma$ denote the set of all letters of $u$. Then the sequence $w \setminus \sigma$ consists of the word $v$. Since $v$ is an assembly word, its length is even, which contradicts Proposition \ref{Prop_equivalent_cond_F2n+1}, Item (4).
\end{Proof}

\section{ Proof of Conjecture \ref{conj_upper_bound_only_for_TC}}
\label{sec:main-result}

In this section, we prove Conjecture \ref{conj_upper_bound_only_for_TC}, which specifies the structure of a maximal assembly word. 
We start showing that each maximal word in some way contains a tangled cord. Then we show that, in order to satisfy Proposition \ref{Prop_equivalent_cond_F2n+1}, a maximal word contains nothing but a tangled cord. We underline that everything is done in terms of assembly words, not assembly graphs, and we only return to the graph terms in Theorem \ref{Th_no_even_is_TCn} and Corollary~\ref{Col_maximal_graph}.

The following lemma shows that a maximal word $w$ always contains a tangled cord that starts with the first letter of $w$ and ends with its last letter.
\begin{Lem}\label{Lem_exists_framing_Tcn}
    Let $w$  be an assembly word in the alphabet $\Sigma$,  $|w|=2n$, and suppose that $w$  cannot be represented as a composition of two assembly words. 
    Then there exist $s \geqslant 1$ and $\sigma = \{t_1, t_2,  \ldots, t_{s}\} \subseteq \Sigma$ such that 
    $\occ_1(t_1) = 1 $, $\occ_2(t_s) = 2n$
    and $w(\sigma) = w_{TC_s}$ up to a renaming of the alphabet letters, i.\,e., $w(\sigma)$ equals $t_1t_2t_1t_3t_2\dots t_s t_{s-1}t_s$. 
\end{Lem}
\begin{Proof}
	We 
    construct such a tangled cord step by step, making the next step while the last letter of obtained tangled cord is not the last letter of~$w$.
	
	We start with the first letter of $w$, set $t_1 = w[1]$.  
	
	If $\occ_2(t_1)=2n$, we set $s = 1$ and the required word is $w_{TC_1} = t_1t_1$.  

 Suppose that we already constructed the subword $$w_{TC_{k}} = w(\{t_1, t_2, \dots, t_{k-1}, t_{k}\})$$ of length $2k$ such that for every $1 \leqslant i \leqslant k$, the letter $t_i$ is chosen in such a way that 
     $$w(\{t_1, t_2, \dots, t_{i-1}, t_{i}\}) = w_{TC_{i}} = t_1t_2t_1t_3t_2\dots t_{i} t_{i-1}t_{i},$$ and $ \occ_2(t_{i})$ is the maximal possible with this property.

If $\occ_2(t_{k})=2n$ then the required tangled cord is constructed. Otherwise we add the next letter as follows. 

	Consider the subword $w[(\occ_2(t_{k})+1) : 2n]$.  According to the   conditions of the lemma, $w$ cannot be a composition of two assembly words. Then there exists a letter $t_{k+1}$ such that $\occ_1(t_{k+1}) < \occ_2(t_k) < \occ_2(t_{k+1})$. 
	
    Note that $1 < \occ_2(t_1) < \occ_2(t_2) < \ldots < \occ_2(t_k) < \occ_2(t_{k+1}).$ There exists $i$, $1 \leqslant i \leqslant k$, such that $\occ_2(t_{i-1}) < \occ_1(t_{k+1}) < \occ_2(t_i)$, here we formally set $\occ_2(t_{0}) = 1$.  Then let us consider $$w(\{t_1, t_2, \dots, t_i, t_{k+1}\}) = t_1t_2t_1t_3t_2\dots t_i t_{i-1}t_{k+1} t_i t_{k+1} = w_{TC_{i+1}}.$$
	Assume firstly that $i < k$. Then up to renaming of letters the words $w(\{t_1, t_2, \dots, t_i, t_{i+1}\})$ and $w(\{t_1, t_2, \dots, t_i, t_{k+1}\})$ are both   equal to $w_{TC_{i+1}}$.  However, $\occ_2(t_{k+1}) > \occ_2(t_{i+1})$, contradicting the choice of $t_{i+1}$. This contradiction implies that $i = k$, so $w(\{t_1, t_2, \dots, t_k, t_{k+1}\}) = w_{TC_{k+1}}$, proving that there exists at least one letter $t_{k+1}$ with the desired property. Now we choose $t_{k+1}$ in such a way that $w(\{t_1, t_2, \dots, t_k, t_{k+1}\}) = w_{TC_{k+1}}$ and $\occ_2(t_{k+1})$ is maximal possible. Thus we are able to continue the process.
	
	Since $w$ has finite length, the described process stops at step $s$, returning $\sigma = \{t_1, t_2, \dots, t_s\}$ such that $\occ_1(t_1) = 1 $, $\occ_2(t_s) = 2n$
	and $w(\sigma) = t_1t_2t_1t_3t_2\dots t_s t_{s-1}t_s =  w_{TC_s}$. This completes the proof.
\end{Proof}

\begin{Def}
    Let $w$ be a word (not necessarily an assembly word) in the alphabet $\Sigma$, $|w| = N$. We say that $w$ contains a framing tangled cord of order $s$ if there exist $t_1, t_2, \dots, t_s \in \Sigma$ such that $w[1] = t_1$, $w[N] = t_s$, and $w(\{t_1, t_2, \dots, t_s\}) = w_{TC_s}$ (up to a renaming of the alphabet letters), i.\,e., it equals $t_1t_2t_1t_3t_2\dots t_s t_{s-1}t_s$. 
\end{Def}

\begin{Note}
    Note that if $w$ contains a framing tangled cord, then it may contain two or more different framing tangled cords by choices of different $s$ and $t_2,\dots,t_{s-1}$. For example, the word $w = 1 2 3 4 1 5 2 6 4 5 3 6$ has framing tangled cords $w(\{1, 3, 6 \}) = 131636$, $w(\{1,4,6\}) = 141646$, $w(\{1,2,5,6\}) = 12152656$.
\end{Note}

\begin{Col}
    An assembly word contains a framing tangled cord if and only if this word cannot be represented as a composition of two assembly words.
\end{Col}
\begin{Proof}
    {\it Necessity}. Follows from Lemma \ref{Lem_exists_framing_Tcn}: an assembly word that cannot be represented as a composition of two assembly words always contains a framing tangled cord.
    
    {\it Sufficiency}. Suppose that an assembly word $w$ contains a framing tangled cord $TC_s$ of order $s$ with letters $t_1, \dots, t_s$. Assume by the contradiction that $w$ can be represented as a composition of two words $w = u\circ v$. Consider $ K_{max} = \max \{ \occ_2(t_i) \, | \, 1 \leqslant i \leqslant s, \linebreak \occ_2(t_i) \leqslant |u|\}$, i.\,e.,  the last number of the position corresponding to a second occurrence of a letter from the framing tangled cord in $u$. Note that $w[K_{max}] = t_k $ for some $k = 1, 2, \dots, s$. We have $k < s,$ since for the framing tangled cord, $\occ_2(t_s)= |w| > |u|.$ Now, for the letter $t_{k+1}$ it follows from the tangled cord structure  that \linebreak $\occ_1 (t_{k+1}) < \occ_2 (t_k) = K_{max} \leqslant |u| < \occ_2 (t_{k+1})$. 
    This contradicts the representation of $w$ as the composition $w=u \circ v$, since the letter $t_{k+1}$ occurs both in $u$ and~$v$. 
\end{Proof}

Below we show that if an assembly word $w$ contains a framing tangled cord, then one can remove some letters of the tangled cord from $w$ so that all the remaining subwords are of even length.

\begin{Ex}
Consider the assembly word $w$ containing a framing tangled cord of order $3$. Assume that the parity of the lengths of the  subwords  between the letters of the  cord are given by Figure~\ref{pic_framing_tcn}. To obtain subwords of even length only, one should remove the letters $t_2$ and~$t_3$.

\begin{figure}[h]
	\begin{center}
		\includegraphics[scale=0.65]{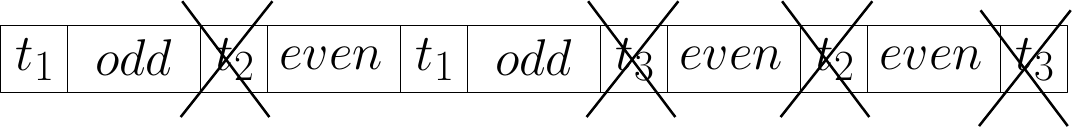}
		\caption{
			Deleting $t_2$ and $t_3$ yields subwords of even length only.}\label{pic_framing_tcn}
	\end{center}
\end{figure}

\end{Ex}

\begin{Lem}\label{Lem_even_subwords_for_framing_TCn}
    Let $w$ be a word of even length (not necessarily an assembly word) in which each symbol occurs at most twice. Suppose $|w| > 2s$ and $w$ contains a framing tangled cord of order $s$, whose letters are $t_1, t_2, \dots, t_s$. Then there exists non-empty $\sigma \subseteq \{t_1, t_2, \dots, t_s\}$ such that all the words in $w \setminus \sigma$ have even length. 
\end{Lem}

\begin{Proof}
The proof is by induction on $s$. Let $|w| = 2n$. 

    {\it The base}: $s = 1$. Hence $w[1] = w[2n] = t_1$, and we set $\sigma = \{t_1\}$. Then $w \setminus \sigma$ contains a unique word, and its length is $2n-2$, which is even.

    {\it Induction step}. Suppose that the lemma's statement holds for all positive integers less than $s$. Then our goal is to prove it for $s$.
    Recall that by the conditions, we have: $w[1] =~t_1$, $w[2n] = t_s$, and $w(\{t_1, t_2, \dots, t_s\}) = w_{TC_s}=t_1t_2t_1t_3t_2\dots t_s t_{s-1}t_s$. We distinguish the following cases.

    {\it Case 1.} $\occ_2(t_i)$ is even for some $1 \leqslant i < s$. Note that $\occ_1(t_i)$ and $\occ_2(t_i)$ are well defined since the framing tangled cord is an assembly word.  Then consider the word $w' = w[1:~\occ_2(t_i)]$ of even length. It contains a framing tangled cord of order $i$. Indeed, $w'(\{t_1, t_2, \dots, t_i\}) = w_{TC_i}=t_1t_2t_1t_3t_2\dots t_i t_{i-1}t_i$. Also, $\occ_1(t_{i+1}) < \occ_2(t_i)$, hence $|w'| > 2i$. Using induction hypotheses, one can find $\sigma \subseteq \{t_1, t_2, \dots, t_i\}$ such that  all the words in $w' \setminus \sigma$ have even length. Then clearly all the words in $w \setminus \sigma$ have even length (since $\occ_2(t_i)$ and $|w|$ are both even), as required, see Figure~\ref{pic_case_occ2_even}.
    \begin{figure}[h]
	\begin{center}
		\includegraphics[scale=0.63]{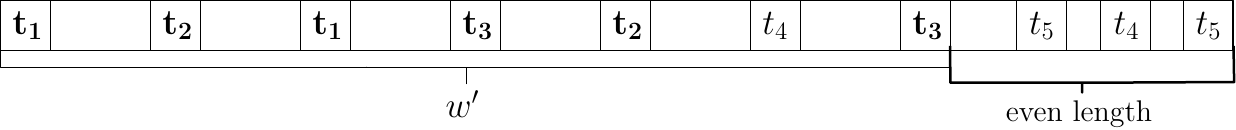}
		\caption{
			Illustration for Lemma~\ref{Lem_even_subwords_for_framing_TCn}, Case 1, $s=5$, $i=3$. The framing tangled cord for $w'$ is bolded.}\label{pic_case_occ2_even}
	\end{center}
\end{figure}

    {\it Case 2.} $\occ_1(t_j)$ is odd for some $1 < j \leqslant s$. This reduces to Case 1 via reversing the word $w$ (all odd indices become even and vice versa, since $|w|$ is even and all first occurrences of symbols become second occurrences). 
    Note that a reversed tangled cord is also a tangled cord, i.\,e., the symmetric argument is valid, see Remark~\ref{note_tc_n_is_symmetric} and Figure~\ref{pic_case_occ1_odd}.

    \begin{figure}[h]
	\begin{center}
		\includegraphics[scale=0.65]{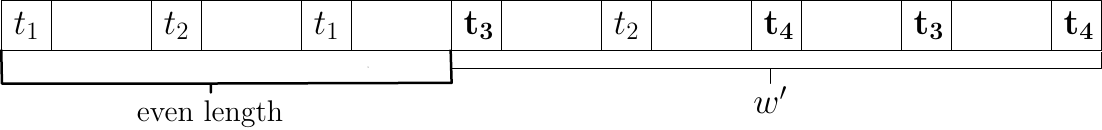}
		\caption{
			Illustration for Lemma \ref{Lem_even_subwords_for_framing_TCn}, Case 2, $s=4$, $j=3$. The framing tangled cord for $w'$ is bolded.}\label{pic_case_occ1_odd}
	\end{center}
\end{figure}

    {\it Case 3.} The remaining case corresponds to the following conditions. The number  $\occ_1(t_j)$ is even for all $1 < j \leqslant s$ and the number $\occ_2(t_i)$ is odd for all $1 \leqslant i < s$. We set $\sigma = \{t_1, t_2, \dots, t_s\}$ and prove that $w \setminus \sigma$ contains only the words of even length, see Figure~\ref{pic_case_occ1_odd1}.

 We note that each of the elements of $w \setminus \sigma$, except the first one \linebreak $w[\occ_1(t_1) +1 : \occ_1(t_2) - 1]$ and the last one $w[\occ_2(t_{s-1}) +1 : \occ_2(t_s) - 1]$, is situated exactly between first and second occurrences of two distinct letters of $TC_s$: it is either \linebreak $w[\occ_1(t_j) +1 : \occ_2(t_{j-1}) - 1]$ for some $1 < j \leqslant s$ or $w[\occ_2(t_{i}) +1: \occ_1(t_{i+2}) - 1]$ for some $1 \leqslant i < s - 1$.
    
    Since for $1< j \leqslant s$ every first occurrence of $t_j$ is even and for $1\leqslant i < s$  every second occurrence of $t_i$ is odd, we have that the number of letters between these occurrences is even. Indeed,
\begin{align*}
    &\left|w[\occ_1(t_j) +1 : \occ_2(t_{j-1}) - 1]\right| = \occ_2(t_{j-1}) - \occ_1(t_j) - 1, \quad \text{  is even for $j > 1$};\\
    &\left|w[\occ_2(t_{i}) +1: \occ_1(t_{i+2}) - 1]\right| = \occ_1(t_{i+2}) - \occ_2(t_i) - 1, \quad  \text{   is even for $i < s-1$}.\\
\end{align*}
    Finally we note that   $\occ_1(t_1) = 1$ and $\occ_1(t_2)$ is even, hence  $$\left|w[\occ_1(t_1) +1: \occ_1(t_2) - 1]\right| = \occ_1(t_2) - \occ_1(t_1) - 1 $$ is even. 
    Similarly, $\occ_2(t_s) = 2n$ and $\occ_2(t_{s-1})$ is odd, hence  $$\left|w[\occ_2(t_{s-1}) +1: \occ_2(t_s) - 1]\right| = \occ_2(t_s) - \occ_2(t_{s-1}) - 1$$ is even. 

    Thus all the words in $w \setminus \sigma$ have even length. This completes the proof.
     \begin{figure}[h]
	\begin{center}
		\includegraphics[scale=0.65]{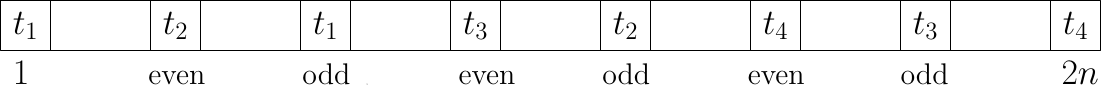}
		\caption{
			Illustration for Lemma \ref{Lem_even_subwords_for_framing_TCn}, Case 3, $s=4$. Odd and even positions alternate, so all the subwords between the cord's letters are of even length. 
   }
   \label{pic_case_occ1_odd1}
	\end{center}
\end{figure}
\end{Proof}

Now we provide a characterization of the maximal words, see Definition~\ref{def_max_word}, and state that a double occurrence word $w$ is maximal if and only if $w=w_{TC_n}$.

\begin{Th}\label{Th_no_even_is_TCn}
    Let $w$ be an assembly word in the alphabet $\Sigma = \{1,2,\dots,n\}$ and $|w| = 2n$. Then the following statements are equivalent:
    
    \begin{itemize}
        \item [1)] For any proper non-empty subset $\sigma \subset \Sigma$, at least one word from $w\setminus \sigma$ has odd length.
        \item [2)]   $w = w_{TC_n}$ up to renaming of symbols. 
    \end{itemize}
\end{Th}
\begin{Proof}
    We prove the equivalence consecutively applying facts from above.
    \begin{itemize}
	 	\item [$1 \Rightarrow 2$]
            First, we note that according to Proposition~\ref{Prop_equivalent_cond_F2n+1}, $w$ is a maximal word (see Definition \ref{def_max_word}). Hence, due to Corollary~\ref{Col_max_not_composition}, $w$ cannot be represented as a composition of two assembly words. From Lemma~\ref{Lem_exists_framing_Tcn} it follows that  $w$ contains a framing tangled cord $TC_s$ of order $s$. Finally, if $|w| > 2s$, then Lemma~\ref{Lem_even_subwords_for_framing_TCn} yields that $w$ is not maximal, which contradicts Item 1. Hence $|w| = 2s = 2n$, thus $w = w_{TC_n}$ up to renaming of symbols. 
            
            \item [$2\Rightarrow 1$]
            It is known (see \cite[Corollary 4.4]{Four_Reg_Graphs_w_Rigid_Vert_DNA_Rec_2013}) that $w_{TC_n}$ is a maximal word. Hence Item 1 follows from Proposition~\ref{Prop_equivalent_cond_F2n+1}.
    \end{itemize}
\end{Proof}

Application of Theorem~\ref{Th_no_even_is_TCn} to the result of  Proposition~\ref{Prop_equivalent_cond_F2n+1} implies our main result. 

\begin{Col} \label{Col_maximal_graph}
    Let $\Gamma$ be a simple assembly graph, $\mathcal{V}(\Gamma) = n$. 
    The number of all Hamiltonian sets of polygonal paths for the graph $\Gamma$ is maximal, i.e., $|\mathcal{C}(\Gamma)| = F_{2n+1}-1$ if and only if   $\Gamma = TC_n$ (up to an isomorphism).
\end{Col}

This proves Conjecture \ref{conj_upper_bound_only_for_TC}. As a corollary, we obtain that each non-maximal word can be split into subwords of even length via removal of some tangled cord. In fact, the corollary below states that the subword formed by the minimal set of symbols that can be removed from a word such that the remaining subwords are of even length always forms a tangled cord.

\begin{Col} \label{Col_we_removed_maximal_word}
    Let $w$ be an assembly word in the alphabet $\Sigma$, $|w| = 2n$. Suppose that $w \neq w_{TC_n}$ (so it is not a maximal word), and let $\sigma \subset \Sigma$ be a non-empty subset of letters of minimal size such that $w \setminus \sigma$ is a collection of subwords of even length only. Then $\Gamma_{w(\sigma)}$ is a tangled cord.
\end{Col}
\begin{Proof}
    Suppose that $|\sigma| = k$ and $w(\sigma) = s_1 s_2 \dots s_{2k}$
    is not a tangled cord. Note that each letter from $\sigma$ occurs in $w(\sigma)$ exactly twice, i.e., $s_i$ are not all distinct. 
    Consider the following representation of $w$:
    $$w = S_1\,s_1\,S_2\,s_2 \dots S_{2k} \, s_{2k} \, S_{2k+1}.$$
    Here for each $i\in \{1,2, \dots, 2k+1\}$ we have that $S_i$ is either an empty subword or $S_i\in w \setminus \sigma.$

    Note that $w(\sigma)$ is not a maximal word by Theorem \ref{Th_no_even_is_TCn}. Then, according to Proposition~\ref{Prop_equivalent_cond_F2n+1}, there exists a proper non-empty subset $\gamma \subset \sigma$ such that $w(\sigma) \setminus \gamma$ consists of subwords of even length only.

    Let $|\gamma| = l < k$ and note that $w(\gamma)$ represents as $w(\gamma) = (w(\sigma))(\gamma) = g_1 g_2 \dots g_{2l}$, $g_i\in \Sigma$. 
    We have a representation of $w(\sigma):$
    $$w(\sigma) = G_1 g_1 G_2 g_2 \dots G_{2l} g_{2l}  G_{2l+1}.$$
    Here for each $j\in \{1,2, \dots, 2l+1\}$ we have that $G_j$ is either an empty subword or $G_j\in w (\sigma) \setminus \gamma.$ Note that $G_1 = s_1 s_2\dots s_{|G_1|}$, $g_1 = s_{|G_1|+1}$, in general, $g_i = s_{|G_1|+\dots + |G_i|+i}$ and $$G_i = w(\sigma)\left[|G_1|+\dots + |G_{i-1}|+ (i-1)+1: |G_1|+\dots + |G_i|+i -1 \right].$$  

    Finally, consider the representation of the initial word $w$ with the letters from $\gamma$ and subwords:
    $$w = U_1 g_1 U_2 g_2 \dots U_{2l} g_{2l} U_{2l+1}.$$
    Here for each $t\in \{1,2, \dots, 2l+1\}$ we have that $U_t$ is either an empty subword or $U_t\in w \setminus \gamma.$ E.\,g., $U_1 = S_1s_1S_2s_2\dots s_{|G_1|}S_{|G_1|+1}$. Note that $|U_i|$ is even. To be more precise,
    \[ 
    |U_i| = |G_i| + \sum_{r = |G_1| + \ldots + |G_{i-1}| + i}^{|G_1| + \ldots + |G_{i}| + i} |S_r|, \quad i = 1, 2, \dots, 2l+1.
    \]
    Here each $|S_r|$ is even by the conditions, $|G_i|$ is even by the choice of $\gamma$ and it is  the number of letters from $s_1, s_2, \dots, s_{2k}$ that are presented in $U_i$.
    
    Therefore, $w \setminus \gamma$ consists of subwords of even length only, which contradicts the fact that the size of $\sigma$ is minimal. The obtained contradiction implies that  $\Gamma_{w(\sigma)}$  is the tangled cord $TC_k$. 
\end{Proof}

\section*{Acknowledgments}
The authors are cordially thankful to the referees for valuable suggestions improving the presentation of the results and for the important references. 

The research of E. Kreines was partially supported by ISF Grant  1092/22, also she is  grateful to colleagues at Tel Aviv University for a warm working atmosphere. 
N. Jonoska was partially supported by the grants NSF DMS-2054321, CCF-2107267, CCF-2505771, and the W.M. Keck Foundation.
The work of A.~Maksaev was supported by the HSE University Basic Research Program.


\end{document}